\documentclass[11pt]{article}
\usepackage{amsmath}
\usepackage{amssymb}
\usepackage{amscd}
\usepackage{latexsym}
\usepackage{theorem}
\usepackage{setspace}
\usepackage{epsfig}
\usepackage{amsfonts}
\usepackage{enumerate}
\usepackage{xypic}
\usepackage{mathrsfs}

\newtheorem{theorem}{Theorem}[section]

\newtheorem{proposition}[theorem]{Proposition}
\newtheorem{corollary}[theorem]{Corollary}

{
\theorembodyfont{\rmfamily}
\theoremstyle{plain}
\newtheorem{definition}[theorem]{Definition}
\newtheorem{example}[theorem]{Example}

\newtheorem{remark}[theorem]{Remark}
}

\newcommand{\id}{\operatorname{id}}

\renewcommand{\dim}{\operatorname{dim}}

\newcommand{\End}{\operatorname{End}}
\newcommand{\Rep}{\operatorname{Rep}}
\newcommand{\Hom}{\operatorname{Hom}}
\newcommand{\Coh}{\operatorname{Coh}}
\newcommand{\Ext}{\operatorname{Ext}}

\newcommand{\Od}{\displaystyle\bigoplus}
\newcommand{\Pd}{\displaystyle\prod}

\newcommand{\C}{{\mathbb{C}}}
\renewcommand{\P}{{\mathbb{P}}}
\newcommand{\Z}{{\mathbb{Z}}}
\newcommand{\N}{{\mathbb{N}}}
\newcommand{\Q}{{\mathbf{Q}}}

\renewcommand{\O}{\mathcal{O}}
\newcommand{\F}{\mathcal{F}}

\newcommand{\mcC}{\mathcal{C}}

\newcommand{\otn}{\{1,\ldots,n\}}

\newcommand{\bigmid}{\hs\Big{|}\hs}
\newcommand{\subs}{\subseteq}

\newcommand{\hs}{\hspace{3pt}}
\renewcommand{\ss}{\substack}

\renewcommand{\a}{\alpha}

\renewcommand{\mod}{{/\!\!/\!\!}}

\newcommand{\Spec}{\mathrm{Spec}\,}
\newcommand{\Proj}{\mathrm{Proj}\,}

\newcommand{\Mchi}{M_\chi}
\newcommand{\Vchi}{V_\chi}

\newcommand{\Bchi}{B_\chi}
\newcommand{\Echi}{E_\chi}
\newcommand{\Tchi}{T_\chi}

\newcommand{\pichi}{\pi_\chi}

\newcommand{\chiss}{R_{\a}(\Q)^{\,\chi-ss}}
\newcommand{\chist}{R_{\a}(\Q)^{\,\chi-st}}
\newcommand{\shom}{\mathcal{H}om}
\newcommand{\res}{\tilde{R}_{\a}(\Q)}
\newcommand{\resp}{\tilde{R}_{\a}(\Q')}
\newcommand{\Ov}{\mathcal{O}_V}
\newcommand{\gm}{\mathbb{G}_m}
\newcommand{\bij}{\beta_{ij}}
\newcommand{\Path}{P}
\renewcommand{\k}{\frak k}
\newcommand{\Db}{\mathcal{D}}

\newcommand{\GL}{\operatorname{GL}}

\newcommand{\roq}{\mathcal{R}ep_\alpha{\Q}}
\newcommand{\roqp}{\mathcal{R}ep_\alpha{\Q'}}

\renewcommand{\(}{\left(}
\renewcommand{\)}{\right)}

\newcommand{\qed}{\hfill \mbox{$\Box$}\medskip\newline}
\newenvironment{proof}{\noindent {\bf Proof:}}{\qed \par}

\begin{document}

\noindent {\LARGE \bf Moduli spaces for Bondal quivers}
\bigskip\\
{\bf Aaron Bergman}\footnote{Supported by the National Science Foundation under Grant Nos. PHY-0071512 and PHY-0455649, and the US Navy, Office of Naval Research, Grant Nos. N00014-03-1-0639 and N00014-04-1-0336, Quantum Optics Initiative.}\\
Department of Physics\footnote{Preprint number UTTG--16--05.},
University of Texas,
Austin, TX 78712\smallskip \\
{\bf Nicholas Proudfoot}\footnote{Supported
by a National Science Foundation Postdoctoral Research Fellowship.}\\
Department of Mathematics, University of Texas,
Austin, TX 78712
\bigskip
{\small
\begin{quote}
\noindent {\em Abstract.}
Given a sufficiently nice collection of sheaves on an algebraic variety $V$,
Bondal explained how to build a quiver $Q$ along with an ideal of relations in
the path algebra of $Q$ such that the derived category
of representations of $Q$ subject to these relations is equivalent to the derived category of coherent
sheaves on $V$.  We consider the case in which these sheaves are all
locally free and study the moduli spaces of semistable representations
of our quiver with relations for various stability conditions.  
We show that $V$ can often be recovered as a connected
component of such a moduli space and we describe the line bundle induced by a GIT construction of the moduli space in terms of the input data.  In certain special cases,
we interpret our results in the language of topological string theory.
\end{quote}
}
\bigskip

\noindent
An algebraic variety $V$ is completely determined by the abelian category
$\Coh(V)$ of coherent sheaves on $V$ \cite{Ga}, and it is therefore a natural problem
to find a way to describe this category in concrete terms.  If $V$ is affine, then $\Coh(V)$
is nothing more than the category of finitely generated modules over the 
algebra of global functions on $V$.  
If we have a presentation of this algebra, this may be interpreted
as a `presentation' of the category $\Coh(V)$.
In the projective case, it is unreasonable to expect $\Coh(V)$ to be equivalent
to the category of modules over any ring.  It is sometimes the case, however,
that such an equivalence can be constructed after passing to the % bounded
derived category $\Db\Coh(V)$.  The derived category carries less information
than the abelian category $\Coh(V)$, but it is enough to reconstruct such invariants
as cohomology, K-theory, and higher Chow groups, as well as a great
deal of information about the birational geometry of $V$.
If $V$ is Calabi-Yau, then an object of $\Db\Coh(V)$ may be thought of as
a D-brane in type IIB topological string theory on $V$ \cite{AD,Do,Sh}.  This category
is therefore of significant physical interest, and is a fundamental ingredient
in the formulation of homological mirror symmetry \cite{Ko}.

Let us describe more concretely how one might attempt to construct such an equivalence.
Given an object
$E$ of $\Coh(V)$, there is a natural functor $F$ from $\Coh(V)$ to the category
of finitely generated right modules over the endomorphism algebra
$\End(E)$, or left modules over the opposite algebra $\End(E)^\text{op}$,
taking a sheaf $\mathcal F$ to the module $\Hom(E,\mathcal F)$.
This functor will almost never be either faithful or essentially surjective, but
if $E$ satisfies certain technical conditions,
then Rickard shows that the right derived functor $RF$ from $\Db\Coh(V)$ to the derived category of left modules over $\End(E)^\text{op}$ will be
an equivalence.  (See Definition \ref{compactspan} 
and Theorem \ref{equiv} for more details.)
If $E$ decomposes as a direct sum of smaller objects
$E = \oplus_{i=1}^n E_i$, then $\End(E)^\text{op}$ may be expressed as the path
algebra of a quiver with $n$ nodes, modulo certain relations (which may not be admissible).
One should think of the description of such a quiver 
along with its relations as an analogue of a presentation of the coordinate ring of an affine variety. 

Much work has gone into finding such collections of sheaves on projective
varieties.  The goal of this paper is not to
find these collections, but rather to assume that one is given, and to study
various moduli spaces of representations of the corresponding algebra.
If the sheaves in the collection are vector bundles, 
there is a tautological map from $V$ to the moduli
stack of quiver representations, taking a point $p$
to $F(\O_p)$, a representation in which the vector space associated
to the node $i$ is equal to the dual of the fiber of $E_i$ at $p$.
Thus, we restrict to representations in which the dimension of the vector space at node $i$ is
the rank of the vector bundle $E_i$. % In this case, if $RF$ is an equivalence,
% then the tautological map is an open embedding (Theorem \ref{bz}).
Our goal is to consider coarse moduli {\em spaces} of semistable representations
% of $Q$
for various choices of stability condition and to relate these spaces to $V$.
The representation stack may be presented as the quotient of an affine
variety by the action of an algebraic group $G$, so these moduli spaces
can be constructed as geometric invariant theory (GIT) quotients
with respect to some character $\chi$ of $G$.
In general, $V$ need not map to such a space as the representations in the image
of the tautological map may not be semistable.  Even if $V$ does map
to one of these moduli spaces, the map may not be an inclusion, as
representations whose closures in the stack intersect in a semistable
representation are identified in the moduli space.
In Section \ref{stable}, the main section of this paper, we address the problem
of determining when this map exists, and when it does, we study its structure.
The GIT construction gives us not just a moduli space, but a moduli space equipped
with an ample line bundle.
Under suitable hypotheses, we show that $V$ may be identified with
(a connected component of) the moduli space of stable quiver representations,
and we identify the induced line bundle on $V$ in terms of $\chi$ and the
vector bundles with which we started (Theorem \ref{gn}).

Section \ref{zero} is devoted to a case of physical interest, in which $V$ is the total
space of the dual of an ample line bundle $L$ on a projective variety $X$, and the collection
on $V$ is pulled back from a particularly nice collection of line bundles on the base.
In this case we study the affine quotient $M_0$, and prove that it has an irreducible
component whose canonical reduced subvariety is isomorphic to $V_0$,
the affine variety obtained from $V$ by collapsing the zero section of the bundle.
If $L$ is the anticanonical
bundle on a Fano survace, this result may be interpreted in the language of topological 
string theory as in the physics paper \cite{BP}.
The quiver moduli space $M_0$ parameterizes ground states of a quantum
field theory that describes the behavior of open strings ending on a certain
D-brane supported at the tip of the Calabi-Yau cone $V_0$.
In general, the quantum field theory associated to a D-brane contains
fields which are sections of the normal cone to the support of the D-brane.
In this case, the normal cone is $V_0$ itself, and a section is simply a point in $V_0$.
For physical reasons, it is expected that the space of sections of this normal cone
should be a component of the moduli space of vacua in the quantum field theory,
and therefore that $V_0$ should be a component of the quiver moduli space.
Up to the issue of reducedness of $M_0$, this is now a theorem.

A special case of the situation discussed above occurs when the collection on $X$ 
is a simple helix (see Example \ref{bridge}).  
% In this case, let $V$
% be the total space of the canonical bundle of $X$, and let $V_0$ be the
% anticanonical cone of $X$, which may be obtained from $V$ by collapsing
% the zero section.  Then $V_0$ is a Calabi-Yau cone, and $V$ is a crepant
% resolution of $V_0$.  We therefore recover $V_0$ as a component of the moduli space of vacua
% for a certain supersymmetric Yang-Mills theory.
In Section \ref{resolutions}, we 
% use the results of Section \ref{zero} 
construct the Fano variety $X$ with its anticanonical
line bundle as a GIT quotient of a smooth variety with respect to a canonical polarization
(Theorem \ref{smoothgit}).  In particular, we obtain a result along the lines of those
of Section \ref{stable} while eliminating the dependence on
the choice of character $\chi$.

\paragraph{\bf Acknowledgments.}  The authors are grateful to 
David Ben-Zvi, Brian Conrad, Deepak Khosla, and Gregory Smith
for invaluable discussions.  We are especially grateful to Alastairs Craw and King for
their detailed comments and suggestions.

\begin{section}{Bondal quivers}\label{bondal}
Let $Q$ be a directed graph with finitely many nodes $\otn$, and let $\k$
be an algebraically closed field.
Let $\Path_{ij}(Q)$ denote the $\k$-vector space spanned by the set of all paths
in $Q$ from the node $i$ to the node $j$, including the path of length
zero at each vertex.
The direct sum 
$\Path(Q) = \bigoplus \Path_{ij}(Q)$ is naturally an algebra over $\k$ with multiplication
$\Path_{jk}\otimes \Path_{ij}\to \Path_{ik}$
given by concatenation of paths.
Let $I\subs \Path(Q)$ be a two-sided homogeneous ideal contained in the square
of the ideal of paths of nonzero length; such an ideal is called \textit{admissible}.  
The pair $\Q = (Q,I)$ with $I$ admissible is called a {\em quiver with relations}.  
The algebra $P(\Q) := \Path(Q)/I$ is called the {\em path algebra} of $\Q$,
and inherits a grading $P(\Q) = \bigoplus P_{ij}(\Q).$

To any quiver with relations $\Q$ we may associate a $\k$-linear category $\mcC(\Q)$
with objects $\otn$, and morphisms from $i$ to $j$ equal to $P_{ij}(\Q)$.
% Put differently, $\mcC(\Q)$ is the category with $P(\Q)$ as its ring of morphisms.
A {\em representation} of $\Q$ is defined to be a functor of $\k$-linear
categories from $\mcC(\Q)$ to the category $\operatorname{Vect}_{\k}$
of $\k$-vector spaces.  
Equivalently, it is a left module
over the path algebra $P(\Q)$.
Let $\Rep(\Q)$
denote the abelian category of representations of $\Q$ that are
finitely generated as $P(\Q)$-modules.

Let $\mcC$ be a $\k$-linear abelian category and consider
% such that the $\Hom$ space between any two objects is a finitely generated
% right module over the endomorphism algebra of the first object.
a finite collection $E_1,\ldots,E_n$ of objects in $\mcC$. 
% the vector space $\Ext^1_\mcC(E_i,E_j)$
% is finite dimensional for each $i,j$.
The algebra 
$$A := \End_\mcC(\oplus E_i)^\text{op}$$ has a distinguished collection $\{e_i\}$
of idempotents, where $e_i$ acts as $\delta_{ij}$ times the identity endomorphism on $E_j$.
Suppose that $A$ is equipped with a grading by the natural numbers, 
with each graded piece finite dimensional, and that
the degree zero part $A_0$ is spanned by the idempotents $\{e_i\}$.
It then makes sense to define the one dimensional representation
$$S_i := A\Big{/} A_{+} + \k\{e_j\mid j\neq i\}$$
on which $e_i$ acts as the identity, and all other idempotents and all elements
of positive degree act by zero.
Let $Q$ be a quiver on $n$ vertices
with arrows from $i$ to $j$ given by a basis for
$\Ext^1_A(S_i,S_j)^\vee$.
There is a map
$$\Ext^2_A(\oplus S_i,\oplus S_i)^\vee\to
\bigoplus_{k\geq2}\left(\Ext^1_A(\oplus S_i,\oplus S_i)^\vee\right)^{\otimes k}$$
given by the $A_\infty$ structure on $\Ext^\bullet_A(\oplus S_i,\oplus S_i)$,
whose image generates an admissible ideal $I\subs P(Q)$.  Let $\Q$ be the corresponding
quiver with relations; we refer to $\Q$ as the {\em Bondal quiver} for the collection
$E_1,\ldots,E_n$.  The following proposition may have been known to the experts for some time, 
but the first proof of it of which we are aware has recently been given by Segal \cite[2.13]{Se}.

\begin{proposition}
The path algebra $P(\Q)$ is isomorphic to $A$.
\end{proposition}

There is a natural functor
$$F:\mcC\to\Rep(\Q)$$
taking an object $\F\in\mcC$ to a representation of $\Q$ in which the node $i$
is mapped to the vector space $\Hom_\mcC(E_i,\F)$.  This functor
is left exact, and thus (assuming a nice class of adapted objects) induces a right derived functor
$$RF:\Db(\mcC)\to\Db\Rep(\Q)$$
on (unbounded) derived categories.

\begin{definition}\label{compactspan}
An object $E$ of $\Db(\mcC)$ is called {\em compact} if the functor $\Hom_{\Db(\mcC)}(E,-)$
commutes with infinite direct sums.
If $\mcC = \Coh(V)$ for a algebraic variety $V$ over $\k$, then the compact objects
of $\Db(\mcC)$ are those which are locally quasi-isomorphic to bounded complexes of locally free sheaves.  If $V$ is smooth, this is simply the class of all complexes quasi-isomorphic to a bounded complex. The derived category $\Db(\mcC)$ is said to be {\em spanned} by a set of objects if for all nonzero objects $F$ of $\Db(\mcC)$, there exists an object $E$ in that set 
such that $Hom_{\Db(\mcC)}(E,F)\neq 0$.  
% We say that $\Db(\mcC)$ is {\em compactly generated}
% if it is spanned by the collection of all perfect objects.
\end{definition}

\begin{theorem}\label{equiv}{\em\cite[6.4]{Ri}}
Suppose that the objects $E_1,\ldots,E_n$ are compact objects that span
$\Db(\mcC)$,
% that $\Db(\mcC)$ is compactly generated, 
and that for all $i,j$, we have $\Ext^k_\mcC(E_i,E_j)=0$ for all $k\neq 0$.
Then $RF$ is an equivalence of triangulated categories.
\end{theorem}
 
% \begin{remark}
The case in which we will be interested is that in which
$\mcC = \Coh(V)$ is the category of coherent sheaves on 
a (not necessarily smooth) algebraic variety $V$ over $\k$. 
% This means that $\Db(\Coh(V))$ is compactly generated, and
% In this case, Rickard's theorem also gives us
% % Rickard's theorem further tells us that we have 
% an equivalence of perfect subcategories
% $$\mathcal{D}^b(Coh(V)) \cong \mathcal{D}^b(\Rep(\Q))$$
% between the bounded derived categories.
% \end{remark}
%, and $E_1,\ldots,E_n$ are line bundles.
In order to endow $\End_V(\oplus E_i)^\text{op}$ with an appropriate grading
from which to construct a Bondal quiver,
we need some extra structure on $V$ and extra conditions on the sheaves $E_1,\ldots,E_n$.

\begin{definition}  A variety $V$ equipped with an action of the multiplicative group $\gm$
is called {\em nearly projective} if it is projective over its affinization 
$V_0 = \Spec\Gamma(\mathcal O_V)$, 
the $\gm$ action on $V_0$ has a unique fixed point, and $\gm$ retracts $V_0$ to that fixed point.
Algebraically, this means that we may write $V = \Proj R$ for an $\N\!\times\!\Z$-graded ring $R$
with $R_{0,i} = 0$ for $i<0$ and $R_{0,0}\cong\k$.  
Here the $\N$-grading is used to construct the $\Proj$ and the $\Z$-grading
gives the $\gm$ action on $V$.
\end{definition}

\begin{example}
Any projective variety $V$ is nearly projective with respect to the trivial $\gm$ action.
\end{example}

\begin{example}\label{tot}
Suppose that $X$ is projective with an ample line bundle $L^{-1}$,
and let $V$ be the total space of $L$.  Then $V$ is nearly projective with respect to
the scaling action of $\gm$ along the fibers.
\end{example}

\begin{definition}
Let $V$ be nearly projective, and let $E_1,\ldots,E_n$ be $\gm$-equivariant vector bundles on $V$.
We call this collection {\em decent} if $\End(E_i) \cong \Gamma(\mathcal O_V)$ for all $i$, 
$\gm$ acts on the vector space $\Hom(E_i,E_j)$ with non-negative
weights for all pairs $i,j$, and it acts with positive weights if $j<i$.
\end{definition}

Let $A = \End_V(\oplus E_i)^\text{op}$, and write
$$A = \bigoplus_{\substack{1\leq i,j \leq n\\ r\in\Z}} A_{ij}^r,$$
where $A_{ij}^r$ is the $r$-eigenspace of $\Hom(E_i,E_j)$ with respect to the action of $\gm$.
We define a grading on $A$ by assigning degree $j-i + nr$ to $A_{ij}^r$.
The following proposition says that this grading has all of the properties required to define the
Bondal quiver $\Q$ of the collection $\{E_1,\ldots,E_n\}$.

\begin{proposition}\label{okay}
If $V$ is nearly projective and $E_1,\ldots,E_n$ is decent, then
this grading is non-negative, the graded pieces are finite dimensional, and the degree zero part
is spanned by the idempotents $\{e_1,\ldots,e_n\}$.
\end{proposition}

\begin{proof}
Non-negativity follows immediately from decency of $E_1,\ldots,E_n$.
To establish finite dimensionality of the graded pieces, it is sufficient to show that $A_{ij}^r$
is finite dimensional for all $i,j,r$.  Let $\pi_0:V\to V_0$ be the natural projection.
Then $A_{ij}^r$ is equal to the $r$-eigenspace of sections of the sheaf
$(\pi_0)_*\shom(E_i,E_j)$ on $V_0$.  Let us write $V = \Proj R$, and let $R_0$ be the degree
zero piece with respect to the $\N$-grading.  
Then $V_0 = \Spec R_0$, and the $\gm$ action on $V_0$
induces an $\N$-grading on $R_0$ with degree zero piece $R_{0,0}$ equal to $\k$.
A $\gm$-equivariant coherent sheaf on $V_0$ corresponds 
to a finitely generated graded $R_0$-module,
and $A_{ij}^r$ is canonically isomorphic the degree $r$ part, which must be finite dimensional.

The degree zero part of $A$ is equal to the direct sum $\oplus_i A_{ii}^0$.  Since our collection
is decent, $A_{ii} = \End_V(E_i)^\text{op}$ is the free $R_0$-module of rank one generated in degree zero by a single class, namely $e_i$.
\end{proof}

\vspace{-\baselineskip}
\begin{definition}\label{exceptional}
For any $\k$-linear abelian category $\mcC$, an object $E$ in $\mcC$ is called {\em exceptional}
if $\End_\mcC(E)\cong \k$ and $Ext^k_\mcC(E,E) = 0$ for $k\neq 0$.  A collection
$E_1,\ldots,E_n$ is called {\em exceptional} if each $E_i$ is exceptional and 
$\Ext^\bullet_\mcC(E_i,E_j) = 0$ for all $i>j$.  An exceptional collection is called {\em full}
if it spans $\Db(\mcC)$, and {\em strong} if $\Ext_\mcC^k(E_i,E_j) = 0$ for all $k\neq 0$ and all $i,j$.
\end{definition}

\begin{example}\label{ec}
Let $V$ be an irreducible projective variety equipped with the trivial $\gm$ action,
and let $E_1,\ldots,E_n$ be 
a full, strong, exceptional collection of vector bundles on $V$ equipped with the trivial
$\gm$ action.  Then the collection is decent, 
the Bondal quiver makes sense, and the hypotheses of Theorem \ref{equiv} are satisfied.
Such collections are known to exist on
projective spaces of arbitrary dimension \cite{Be}, 
and on all odd dimensional, smooth, quadric hypersurfaces \cite{Ka}.
They are conjectured to exist on complete
flag varieties of semisimple groups \cite[1.2]{Ku}.
King \cite{K2} shows that they exist on all smooth, Fano, toric surfaces,
and Craw and Smith \cite{CS} extend this result to smooth, Fano, toric 3-folds.
Costa and Mir\'o-Roig \cite{CM}  have found more toric examples in arbitrary dimension.
King \cite[9.3 $\&$ 9.4]{K2} conjectured
that such a full, strong, exceptional collection exists on every smooth, projective toric variety,
and (more generally) on {\em any} variety that may be obtained as a GIT quotient
of a vector space by a linear action of a reductive group, provided that a polarization
is chosen for which the notions of stability and semistability coincide.
Hille and Perling \cite{HP} have recently constructed a toric counterexample
to King's conjecture, but the question how common such collections are is still wide open.
\end{example}

\begin{example}\label{bridge}
Let $E_1,\dots,E_n$ be a full, strong, exceptional collection of vector bundles
on a smooth, projective variety $X$, 
and let $L$ be a line bundle on $X$ such that $L^{-1}$ is ample. 
We extend our collection infinitely in both directions
via the formula
\begin{equation*}
E_{i-n} = E_i \otimes L \quad \text{ for all } i\in\mathbb{Z}\,.
\end{equation*}
Such an infinite collection will be called a {\em spiral with respect to $L$}.
A spiral will be called {\em simple} if it satisfies the equation
\begin{equation}
\label{simpspiral}
\Ext^k_X(E_i,E_j) = 0 \quad \text{ for all } k\neq 0 \mathrm{\ and\ } i\le j\ .
\end{equation}
A simple spiral with respect to the canonical bundle with 
$n = \dim X + 1$
will be called a {\em simple helix}\,;
this notion will become important in Remark \ref{gbgs}
and Section \ref{resolutions}.\footnote{\nobreak
Our definition of a simple helix agrees with that of \cite[\S 3]{Br};
the same structure is called a {\em geometric helix}
in \cite[\S 1]{BoP}.  The definition of a helix is inconsistent
in the literature, and we will never use it.
The word {\em spiral} is our own.}

Suppose that $E_1,\ldots, E_n$ generate a simple spiral
with respect to $L$.
Let $V$ denote the total space of $L$, and let
$\pi:V\to X$ be the projection; $V$ is nearly projective by Example \ref{tot}.
For any pair $i,j$, we have
\begin{equation*}
\Hom_V\!\(\pi^*E_i,\pi^*E_j\) 
\cong \Hom_X\!\(E_i,\pi_*\pi^*E_j\)
\cong \Od_{r\geq 0}\Hom_X\!\(E_i,E_j\otimes L^{-r}\),
\end{equation*}
where $r$ is the eigenvalue for the action of $\gm$.  In particular, the collection
$\pi^*E_1,\ldots,\pi^*E_n$ is decent and the Bondal quiver is well-defined.

For any $\mathcal F\in\Db(\Coh(V))$, we have
$$\Ext_V^\bullet(\oplus\pi^*E_i,\mathcal F) =
\Ext_X^\bullet(\oplus E_i,\pi_*\mathcal F),$$
which is trivial if and only if $\mathcal F = 0$.
Hence $\pi^*E_1,\ldots,\pi^*E_n$ span $\Db(\Coh(V))$.
The condition \eqref{simpspiral} ensures that
the bundles $\pi^*E_1,\ldots,\pi^*E_n$
have no higher Ext groups between them,
so Theorem \ref{equiv} applies to this collection.
\end{example}

% \begin{example}\label{stupid}
% Suppose that $V$ is an affine variety and our collection
% consists only of the trivial bundle.  Then $\Q$ has a single node,
% and its path algebra is isomorphic to the coordinate ring of $V$.
% In this case, we do not need to pass to the derived category;
% $F$ itself is already an equivalence.  While this equivalence is very
% important in the study of affine algebraic geometry, it is not at all
% illuminated by the language of quivers.  As a general rule,
% the closer $V$ is to being projective, the less internal structure 
% is carried by the vertices of $\Q$,
% and the more interesting Theorem \ref{equiv} becomes.
% \end{example}

For the rest of the paper we assume that $V$ is nearly projective and that $E_1,\ldots,E_n$
is a decent collection of vector bundles with Bondal quiver $\Q$.  We do {\em not} assume
that the derived functor $RF$ is an equivalence unless we say so explicitly. Let $\alpha_i = \mathrm{rank\ }E_i$ be a vector of natural numbers, and set $\mathcal{R}ep_\alpha(\Q)$ denote the substack of the moduli stack of representations
of $\Q$ for which the vector space associated to the node $i$ has dimension $\alpha_i$.
Over each point in the variety $V$, the fiber of the vector bundle $\bigoplus E_i^\vee$ is naturally a left-module over the algebra $\End_V(\bigoplus E_i)^\mathrm{op} \cong P(\Q)$. Thus, $V$ parametrizes a family of representations, and we have a tautological map:
$$T:V\to\mathcal{R}ep_\alpha(\Q)\ .$$
On the level of points, we have $T(p) = F(\mathcal{O}_p)$,
where $\mathcal{O}_p$ is the structure sheaf of the point $p\in V$.
% Concretely, $T(p)$ is the representation
% of $\Q$ that takes the $i^\text{th}$ node to the vector space 
% $E_i^\vee\!\!\mid_p$, the dual of the fiber of $E_i$ at $p$.
%The {\em tautological map}
%$$T:V\to\mathcal{R}ep_1(\Q)$$
%is defined by associating to any morphism $f:S\to V$ the $\mathcal O_S\otimes\End(\oplus E_i)^{\text{op}}$
%module $\oplus f^*E_i^\vee$.

\begin{theorem}\label{bz}
If $V$ is smooth and $RF$ is an equivalence of derived categories, then $T$ is injective 
and induces an isomorphism on tangent spaces.
\end{theorem}

\begin{proof}
The fact that each $E_i$ is a vector bundle and $\mathcal{O}_p$
has zero dimensional support tells us that all of the higher
right derived functors of $F$ vanish on $\mathcal{O}_p$.
Hence $RF(\mathcal{O}_p) = F(\mathcal{O}_p)$ is an honest
representation, rather than a complex of representations.
% and $T$ gives rise to a family
% of objects in $Rep(\Q)$.
Injectivity of $T$ then follows from the fact that the objects
$\{\mathcal{O}_p\mid p\in V\}$ are all nonisomorphic in $\Db\Coh(V)\cong \Db\Rep(\Q)$ and, hence, in the full subcategory $\Rep(\Q)$.

To see that $T$ induces an isomorphism on tangent spaces, we note that
we have a sequence of isomorphisms
\begin{equation}\label{isoms}
\begin{split}
T_pV &\cong
\Ext^1_V(\mathcal{O}_p,\mathcal{O}_p)\\ 
&\cong \Hom_{\Db\Coh(V)}(\mathcal{O}_p,\mathcal{O}_p[1])\\
&\cong \Hom_{\Db\Rep(\Q)}(RF(\mathcal{O}_p),RF(\mathcal{O}_p)[1]) \\&\cong \Ext^1_{\Rep{\Q}}(F(\mathcal{O}_p),F(\mathcal{O}_p))\\
&\cong T_{T(p)}\roq.
\end{split}
\end{equation}
Let $\mathbb{D} = \mathrm{Spec\ } \k[\epsilon]/\langle\epsilon^2\rangle$.
Then a tangent vector to $V$ is represented by a map $e:\mathbb{D}\to V$, 
a tangent vector to $\roq$ is represented by a family 
of quiver representations over $\mathbb{D}$,
and the differential of $T$ sends $e\in T_pV$ to $\oplus\, e^* E_i^\vee\in T_{T(p)}\roq$.
It remains only to show that this map coincides with the isomorphism of Equation \eqref{isoms}.

Consider the following Cartesian square:
\begin{figure}[h]
\centerline{
\xymatrix{ 
\mathbb{D}\ar[rr]^{e} \ar[d]_{e\times \mathrm{id}}&& V \ar[d]^{\Delta}  \\
V \times \mathbb{D} \ar[rr]_{\mathrm{id}\times e} && V \times V,
}}
\end{figure}

\noindent
and let $\pi : V \times \mathbb{D} \to \mathbb{D}$ and $\rho : V \times \mathbb{D} \to V$ 
denote the projections.
An element of $T_pV \cong \Ext^1_V(\mathcal{O}_p,\mathcal{O}_p)$ 
may be regarded as a family
of coherent sheaves on $V$ parameterized by $\mathbb{D}$, or, equivalently, 
as a coherent sheaf on $V\times \mathbb{D}$.
In these terms, the element represented by $e:\mathbb{D}\to V$
may be identified with the coherent sheaf
$$(\mathrm{id} \times e)^*\Delta_*\mathcal{O}_V 
= (e \times \mathrm{id})_* e^*\mathcal{O}_V.$$
Then the family of quiver representations obtained by applying $F$ is
\begin{equation*}
\begin{split}
\pi_*\mathcal{H}om_{V\times\mathbb{D}}\left(\oplus\, \rho^*E_i, (e \times \mathrm{id})_* e^*\mathcal{O}_V\right)
&= \pi_*(e\times\mathrm{id})_*\mathcal{H}om_\mathbb{D}\left(\oplus\, e^*E_i, e^*\mathcal{O}_V\right) \\
&= \mathcal{H}om_\mathbb{D}\left(\oplus\, e^*E_i, e^*\mathcal{O}_V\right) \\
&= \oplus\,e^* E_i^\vee,
\end{split}
\end{equation*}
which is precisely the tangent vector to $\roq$ obtained by applying the differential of $T$ to $e$.
Thus the isomorphism of Equation \eqref{isoms} is indeed the one induced by $T$.
\end{proof}
\end{section}

\begin{section}{Semistable representations}\label{stable}
As in Section \ref{bondal}, let $E_1,\ldots,E_n$ be a decent collection of $\gm$-equivariant
vector bundles on a nearly projective algebraic variety $V$ over $\k$, and let 
$\alpha_i = \mathrm{rank\ }E_i$. Let $R_{\a}(\Q)$ be the set of representations of $\Q$ in 
which the node $i$ is mapped to a given vector space of dimension $\alpha_i$ the fixed coordinate vector space
$\k^{\a_i}$. This set has the structure of an affine algebraic variety over $\k$. 
% More explicitly,
% consider the natural composition map
%$$\psi_{ijk}:\Hom(E_i,E_j)\otimes\Hom(E_j,E_k)\to\Hom(E_i,E_k),$$
%which induces a dual map
%$$\psi_{ijk}^\vee:\Hom(E_i,E_k)^\vee\to\Hom(E_i,E_j)^\vee\otimes\Hom(E_j,E_j)^\vee.$$
%Then we have 
%\begin{equation}\label{rq}
%R_{\a}(\Q) = \Big\{(a_{ij})\in\Pd_{i,j}\Hom(E_i,E_j)^\vee\bigmid
%\ss{\psi_{ijk}^\vee(a_{ik}) = a_{ij}\otimes a_{jk}\hs\text{ for all $i,j,k$}
%\\\text{and}\quad a_{ii}\(\operatorname{id}_{E_i}\) = 1}
%\Big\}.
%\end{equation}

\begin{example}\label{p2}
Let $V = \P^2$, and let $E_1 = \O, E_2 = \O(1)$, and $E_3 = \O(2)$.
The following picture represents the category $\mathcal C(\Q)$, with each arrow
labeled by the vector space of morphisms between the corresponding objects.

\begin{figure}[h]
\centerline{
\xymatrix{ 
& 2\ar[dl]_{\Gamma(\O(1))}\\
1&& 3\ar[ll]^{\Gamma(\O(2))}\ar[ul]_{\Gamma(\O(1))}
}}
\end{figure}

\noindent
The quiver itself consists of three arrows from $2$ to $1$ and three arrows from $3$ to $2$,
representing bases for the vector space $\Gamma(\O(1))$.  There are no arrows from $3$ to $1$,
because the multiplication map $$\psi:\Gamma(\O(1)) \otimes \Gamma(\O(1)) \to \Gamma(\O(2))$$ is surjective.
An element of $R_{\a}(\Q)$ consists of a pair of vectors 
$a_{12},a_{23} \in\Gamma(\O(1))^\vee$ such that 
$a_{12}\otimes a_{23}$ lies in the image of $\psi^\vee$.
In concrete terms, this means that $a_{12}$ and $a_{23}$ must be parallel.
\end{example}

Let $G = \prod_i \GL(\alpha_i)\big/\gm^{diag}$.  
This group acts naturally on $R_{\a}(\Q)$ by the formula
$$(g_1,\ldots,g_n)\cdot (a_{ij}) = (g_i a_{ij} g_j^{-1}),$$
and two representations are isomorphic if and only if they lie in the same $G$-orbit.
Any representation of $\Q$ in which all nodes are mapped to vector spaces
of the given dimension is isomorphic
to an element of $R_{\a}(\Q)$; this is just the statement that all finite dimensional vector spaces of a given dimension are isomorphic. It follows that the stack $\mathcal{R}ep_\alpha(\Q)$, 
considered in the previous section, is represented 
by the quotient $\left[R_{\a}(\Q)/G\right]$.

Let $\chi = (\chi_1,\ldots, \chi_n)$ be an ordered $n$-tuple of integers satisfying
$\sum\chi_i\alpha_i=0$. We may interpret $\chi$ as a multiplicative character of the group
$G$ by the formula $g\mapsto \det(g_1)^{\chi_1}\!\!\ldots \det(g_n)^{\chi_n}$.
Let $$\Mchi = R_{\a}(\Q)\mod_\chi G$$ be the semiprojective GIT quotient of $R_{\a}(\Q)$
by $G$ with respect to the character $\chi$.  This quotient has two equivalent
interpretations, which we describe below.\footnote{Geometric invariant theory was originally
developed by Mumford \cite{MFK}, but what we need is summarized in the short survey \cite{Pr}.}

Let $B$ be the affine coordinate ring of $R_{\a}(\Q)$.
%Since $R_{\a}(\Q)$ is a subvariety
%of the vector space $\prod\Hom(E_i,E_j)^\vee$, its coordinate ring
%$B$ is a quotient of the polynomial ring
%$\Sym\bigoplus_{i,j}\Hom(E_i,E_j)$.
%(Even
%if this vector space is infinite dimensional, $R_{\a}(\Q)$ is finite dimensional, and therefore
%$B$ is finitely generated.)
The action of $G$ on $R_{\a}(\Q)$ induces an action
on $B$.  For any character $\theta$ of $G$, let $B(\theta)$ be the $\theta$-eigenspace
of $B$, and let $$\Bchi = \Od_{r\geq 0}B(r\chi).$$  The GIT quotient
$\Mchi$ is defined as $\Proj B_\chi$.  This definition makes it clear that
$\Mchi$ is a variety equipped with an ample line bundle, making $\Mchi$
projective over its affinization $M_0 = \Spec B^G$.

An element $a$ of $R_{\a}(\Q)$ is called {\em $\chi$-semistable} if there exists a function
$f\in B(r\chi)$ for some $r>0$ such that $f(a) \neq 0$.  The locus of semistable
points is an open subset of $R_{\a}(\Q)$, and will be denoted $\chiss$.
Such a representation is called
{\em $\chi$-stable} if its stabilizer is finite and its $G$-orbit is closed in $\chiss$.
The locus of stable points is an open subset of $\chiss$ and will be denoted
$\chist$.  Two semistable representations are called {\em S-equivalent}
if the closures of their $G$-orbits intersect in $\chiss$.  There is a surjective
map from $\chiss$ to $\Mchi$ whose fibers are precisely the S-equivalence classes,
so $\Mchi$ may be thought of as the moduli space of semistable representations
of $\Q$ with dimension vector $\alpha$, up to S-equivalence.

Recall the tautological map $T:V\to \roq$.  The variety $\Mchi$ is a quotient
of an open substack of $\roq$, so $T$ induces a rational map
$\Tchi:V\to\Mchi$.
If $\Tchi$ is in fact regular, meaning that 
every tautological representation $T(p)$ is 
$\chi$-semistable, we will say that the character $\chi$ {\em good}.
If in fact $T(p)$ is $\chi$-stable for all $p$, we will say that $\chi$ is {\em great}.

As a first step to analyzing the map $T_\chi$ for various values of $\chi$, we must consider the
case where $\chi = 0$.
In this case, $T\chi$ factors through the affininzation map 
$$\pi_0:V\to V_0 := \Spec\Gamma(\mathcal O_V)$$
via the map $$\varphi_0:V_0\to M_0 = \Spec B^G = \Spec\Gamma(\mathcal O_{\roq})$$ 
obtained by pulling back global functions
from $\roq$ to $V$.
We note that
every element of $R_{\a}(\Q)$ is semistable with respect to the trivial character, so $\chi = 0$
is always good. 

\begin{proposition}\label{first}
The map $\varphi_0:V_0\to M_0$ is a closed embedding.
\end{proposition}

\begin{proof}
This is equivalent to the statement that $T^*:\Gamma(\mathcal O_V)\to\Gamma(\mathcal O_{\roq})$
is surjective.
Choose any node $i$.  The isomorphism $\Gamma(\Ov)\cong\End(E_i)$
coming from decency of the collection allows us to identify the ring of global functions on $V$
with the algebra of loops in $\Q$ based at $i$.  For any function $f\in\Gamma(\Ov)$, let
$s_i(f)$ be the $G$-invariant function on $R_{\a}(\Q)$ taking a representation to $\frac{1}{\a_i}$ times
the trace of the endomorphism obtained by going around the loop corresponding to $f$.
Then $T^*s_i(f) = f$.
\end{proof}

\vspace{-\baselineskip}
\begin{remark}\label{ifone}
If $\a_i=1$, then $s_i$ is a homomorphism, and induces a map $\sigma_i:M_0\to V_0$
of which $\varphi_0$ is a section.  In general, however, $s_i$ fails to be an isomorphism
because trace is not multiplicative.
\end{remark}

For any character $\chi$, consider the line bundle
$$\Echi = \det(E_1)^{\otimes\chi_1}\otimes\ldots\otimes \det(E_n)^{\otimes\chi_n}.$$

\begin{theorem}\label{gn}
Suppose that $V$ is smooth, $RF$ induces an equivalence of derived categories, 
and $\chi$ is great.  
Then $\Tchi$ identifies $V$ with
a connected
component of $\Mchi$, and $\Echi$ with the line bundle induced by the GIT construction.  
\end{theorem}

\begin{proof}
Since $\chi$ is great, $\Tchi$ maps $V$ to the stable locus of $\Mchi$,
which is isomorphic to an open substack of $\roq$.
Theorem \ref{bz} tells us that $\Tchi$ is injective on points and induces an isomorphism on tangent spaces.  Since $V$ is smooth, this implies that $\Tchi$ is an isomorphism onto a
Zariski open subset of $\Mchi$. 

Since $V$ is nearly projective, it is projective over its affinization $V_0$, 
and $\Mchi$ is projective over $M_0$.
Since $\Tchi:V\to\Mchi$ covers the closed immersion $T_0:V_0\to M_0$, its image must be closed.  Thus $\Tchi$ is an isomorphism onto a connected component of $\Mchi$.

To prove the final statement, we note that the character $\chi$ defines an equivariant
structure on the trivial line bundle on $R_\a(\Q)$, which descends to a nontrivial
line bundle $L_\chi$ on the stack quotient $[R_\a(\Q)/G]\cong\roq$.  The GIT line bundle on $\Mchi$
is obtained by restricting $L_\chi$ from $\roq$, so it will suffice to show that $T^*L_\chi = E_\chi$.

Let $\tilde{\mathcal{G}}$ 
be the principal $\prod_i\GL(\a_i)$-bundle on $V$ associated
to the vector bundle $E = \oplus_i E_i$, and let $\mathcal G$ be the principle $G$-bundle
obtained by dividing $\tilde{\mathcal{G}}$ by $\gm^{diag}$.  Then
%, by the definition of the tautological map, 
we have a pullback diagram
of principal $G$-bundles\eject
\begin{figure}[h]
\centerline{
\xymatrix{ 
& \mathcal G\ar[d]\ar[r] & R_\a(\Q)\ar[d]\\
& V\ar[r]^{\!\!\!\!\!\!\!\!\!\!\!\!\! T} & [R_\a(\Q)/G]\, ,
}}
\end{figure}

\noindent
and the line bundles $\Echi$ and $L_\chi$ are the line bundles associated to these principle
bundles via the one-dimensional representation of $G$ given by the character $\chi$.
The statement follows.
\end{proof}

\vspace{-\baselineskip}
\begin{remark}
More generally, the rational map $\Tchi:V\to\Mchi$ factors through the rational map
$\pichi:V\to\Vchi$ via a third rational map $\varphi_\chi:\Vchi\to\Mchi$.
The maps $\pichi$ and $\varphi_\chi$ will both be regular if and only if $\chi$ is good.
\end{remark}

\begin{remark}\label{cs}
Craw and Smith \cite{CS} obtain a result similar to Theorem \ref{gn}, but with 
different hypotheses.  The most
important differences are that they restrict to collections of line bundles and that they assume
that $V$ is toric. In exchange, they are able to substantially weaken the assumption
that $F$ is an equivalence of categories.  
\end{remark}

The remainder of this section will be devoted to giving sufficient criteria for $\chi$
to be good or great in the case where each $E_i$ is a line bundle.
There is a simple description of stability and
semistability of quiver representations due to King \cite[\S 3]{K1}.
For any subset $S\subs\otn$, let 
$\chi_S = \sum_{i\in S}\chi_i$.  We define the {\em support} of a representation
of $\Q$ to be the set of nodes that map to nonzero vector spaces.
A representation $a\in R_{\a}(\Q)$ has a subrepresentation with support $S$
if and only if $a_{ij}=0$ for all $i\in S^c$ and $j\in S$.
King tells us that $a$ is $\chi$-semistable if and only if
$\chi_S \leq 0$ for all supports $S$ of subrepresentations of $a$,
and $a$ is $\chi$-stable if equality is obtained only by the trivial representation
and $a$ itself.

Let $\{m_{ij}\}$ be a collection of non-negative integers,
and define $\chi$ by the formula 
$$\chi = \sum_{i,j}\,m_{ij}\cdot(0,\ldots,0,-1,0,\ldots,0,1,0,\ldots,0),$$
where $-1$ appears in the $i^\text{th}$ spot, and $1$ in the $j^\text{th}$ spot.
Equivalently, we put
$$\chi_\ell = \sum_{i=1}^n m_{i\ell} - \sum_{j=1}^n m_{\ell j}\quad
\text{for all $\ell\leq n$.}$$

\begin{proposition}\label{good}
If $\shom(E_i,E_j)$ is generated by global sections for all $i,j$ such that
$m_{ij}\neq 0$, then $\chi$ is good.
\end{proposition}

\begin{proof}
Let $S$ be the support of a subrepresentation
of $T(p)$ for some $p\in V$.  We need to show that $\chi_S\leq 0$, where
\begin{equation}\label{difference}
\chi_S = \sum_{\ell\in S}\chi_\ell = \sum_{\ell\in S}\sum_{i=1}^nm_{i\ell}
-\sum_{j=1}^n\sum_{\ell\in S} m_{\ell j}.
\end{equation}
The condition that $\shom(E_i,E_j)$ is generated by global sections
says exactly that $T(p)_{ij} \neq 0$ for all $p\in V$.  
Thus if $m_{ij}\neq 0$ and $j\in S$,
$i$ must be in $S$ as well.  
This tells us that every term that appears with a plus
sign above also appears with a minus sign, therefore $\chi_S\leq 0$.
\end{proof}

% \vspace{-\baselineskip}
% \begin{remark}
% This proposition can be understood more conceptually as follows. In order for a character $\chi$ to be good, for every point $p\in V$, there must be an element of $\Bchi$ of nonzero degree that doesn't
% vanish on the quiver representation $T(p)$. Since we are dealing with collections of line bundles, paths immediately give rise to semi-invariants. As each $\shom(E_i,E_j)$ is generated by global sections, there exists a global section which does not vanish at a given point $p$. Since the product of nonvanishing things is nonvanishing, the product of such sections gives a semi-invariant that does not vanish at $p$.
% \end{remark}

We will say that $\{m_{ij}\}$ is {\em sufficient} if the following two conditions are satisfied:
\begin{enumerate}
\item $\shom(E_i,E_j)$ is generated by global sections for all $i,j$ such that
$m_{ij}\neq 0$.
\item It is possible to get from any
one vertex of $Q$ to any other by traveling forward along paths from $j$ to $i$ such
that $\shom(E_i,E_j)$ is generated by global sections, and backward along paths
from $j$ to $i$ such that $m_{ij}\neq 0$.
\end{enumerate}

\begin{proposition}\label{great}
If $\{m_{ij}\}$ is sufficient, then $\chi$ is great.
\end{proposition}

\begin{proof}
Let $S$ be the support of a nonzero subrepresentation of $T(p)$ for some $p\in V$,
and suppose that $\chi_S = 0$.  We need to show that $S = \otn$.
If $j\in S$ and $\shom(E_i,E_j)$ is generated by global sections, then $T(p)_{ij}\neq 0$,
and therefore $i\in S$.  If $i\in S$ and $m_{ij}\neq 0$,
then $-m_{ij}$ is a summand in Equation \eqref{difference}.  
Since every positive summand is canceled by a negative one and $\chi_S = 0$,
the term $m_{ij}$ must appear as well, hence $j\in S$.  In this manner, we can conclude
that the set $S$ is closed under the two operations described in Condition 2 above.
Since $S$ is nonempty, it must contain all of $\otn$.
\end{proof}

\vspace{-\baselineskip}
\begin{remark}\label{gbgs} 
Suppose that $E_1,\ldots,E_n$ generate a simple helix on a projective variety $X$, 
as in Example \ref{bridge}.
Bondal and Polishchuk \cite[2.5]{BoP} show that for all $1\leq i\leq n-1$,
the object $R_{E_{i+1}}E_i \in\Db\Coh(X)$ defined by the exact
triangle $$E_{i+1}\to\Hom(E_i,E_{i+1})^\vee\otimes E_{i+1}\to R_{E_{i+1}}E_i$$
is pure; in other words, it lies in the abelian subcategory $\Coh(X)$.
This is equivalent to the statement that the first map in the triangle is injective,
or that its dual is surjective.  This in turn is the statement that
$\shom(E_i,E_{i+1})$ is generated by global sections, and therefore so is
$\shom(E_i,E_j)$ for all $i\leq j$.  
Furthermore, they prove that the endomorphism
algebra $\End(\oplus E_i)$ is multiplicatively generated by elements of
the vector spaces $\Hom(E_i,E_{i+1})$, which implies that $1$ is the unique
source of $\Q$, and $n$ is the unique sink.
In this case, therefore, $\{m_{ij}\}$ is sufficient if and only if $m_{1n}>0$.
\end{remark}

Remark \ref{gbgs} gives us many examples of characters $\chi$ that satisfy
the hypotheses of Theorem \ref{gn} when our collection generates a simple helix.
This will apply to the collection $\O,\O(1)\ldots,\O(n)$ on $\P^n$ \cite{Be},
as well as to collections on odd dimensional quadrics \cite{Ka}.
We conclude this section with an example in which Remark \ref{gbgs}
does {\em not} apply, but Propositions \ref{good} and \ref{great} do.

\begin{example}\label{dp1}
Let $V$ be the Hirzebruch surface $\mathbb F_1$, the blow-up of $\P^2$ at a single point.
Consider the collection $E_1 = \O, E_2 = \O(D), E_3 = \O(H)$, and $E_4 = \O(2H)$,
where $H$ is the proper transform of  a hyperplane class in $\P^2$ and $D$ is the exceptional divisor.
This collection is full, strong, and exceptional, and has the following Bondal quiver, where
the integers above the arrows indicate the number of distinct arrows between the two nodes
(unlabeled arrows occur with multiplicity one).  
There are nontrivial relations among paths from $4$ to $1$ and among maps from $4$ to $2$.
\begin{figure}[h]
\centerline{
\xymatrix{ 
& 2\ar[dl]\\
1&& 3\ar[ll]\ar[ul]_2 && 4\ar[ll]_3
}}
\end{figure}

\noindent
The only nonzero path in $\Q$ corresponding to a $\shom$ sheaf which is {\em not} generated by global sections
is the arrow from 2 to 1.  
Let $m_{14} = m_{23} = 1$, and set all other $m_{ij}$ equal to zero.
Then $\{m_{ij}\}$ is sufficient, so $\chi = (-1,-1,1,1)$ is great.  
Then Theorem \ref{gn}
tells us that $\Mchi$ has a connected component 
that is isomorphic to $V$ in its projective embedding 
given by the anticanonical bundle $\Echi = \O(3H-D)$.
\end{example}
\end{section}

\begin{section}{D-branes at the tip of a cone}\label{zero}
In this section we continue to assume that $\a_i=1$ for all $i$.
Suppose that a collection $E_1,\ldots,E_n$ of line bundles
generates a simple spiral with respect to another line bundle $L$
on a smooth projective variety $X$,
as in Example \ref{bridge}.  Let $V$ denote the total space of $L$, and let $\Q$ be
the Bondal quiver for the decent collection $\pi^*E_1,\ldots,\pi^*E_n$ on $V$.
Theorem \ref{affine} and Corollary \ref{component} generalize the main result of \cite{BP}.

\begin{theorem}\label{affine}
The map $\varphi_0$ is generically an isomorphism.  More precisely, there exists
a dense open subset $U\subs V_0$ such that $\varphi_0|_U$ is an isomorphism onto its
image, which is open in $M_0$.
\end{theorem}

\begin{proof}
We first observe that for all $i$ and $j$, there exist elements
$$p_{ij}\in\Hom_V\!\(\pi^*E_i,\pi^*E_j\)\cong \Od_{r\geq 0}\Hom_X\!\(E_i,E_j\otimes L^{-r}\)$$
and $$q_{ij}\in\Hom_V\!\(\pi^*E_j,\pi^*E_i\)\cong \Od_{r\geq 0}\Hom_X\!\(E_j,E_i\otimes L^{-r}\)$$
with nonzero
product $\bij = q_{ij}\cdot p_{ij}\in\End_V(\pi^*E_i)\cong\Gamma(\mathcal{O}_V)$.
This follows from ampleness of $L^{-1}$, which ensures that the vector spaces on
the right will be large for large values of $r$.

Recall from Remark \ref{ifone} that 
for each node $i$ we have a homomorphism
$s_i:\Gamma(\Ov)\to\Gamma(\mathcal O_{M_0})$
inducing a map $\sigma_i:M_0\to V_0$ 
with the property that % $\varphi_0^* s_i(f) = f$ for all $f\in\Gamma(\Ov)$.
$\sigma_i\circ\varphi_0 = \id_V$.
Le Bruyn and Procesi \cite[Thm 1]{LP} show that the images of $s_1,\ldots,s_n$ generate
$\Gamma(\mathcal{O}_{M_0})$.
Furthermore, for any element $r\in\Hom(E_j,E_j)$, we have
$$s_i(\bij)\cdot s_j(r) = s_i(q_{ij}\cdot p_{ij})\cdot s_j(r) = s_i(p_{ij}\cdot r\cdot q_{ij}).$$
This means that $s_i$ becomes surjective after inverting the elements 
$s_j(\bij)\in\Gamma(\mathcal O_{M_0})$ for all $j$.
Geometrically, this tells us that there exists a dense open set $U_i$ of $V_0$ (the set on
which $0\neq\varphi_0^*s_j(\bij)=\bij$ for all $j$) over which $\sigma_i$ is an isomorphism.
Since $\varphi_0$ is a section of $\sigma_i$, we are done.
\end{proof}

\vspace{-\baselineskip}
\begin{corollary}\label{component}
The map $\varphi_0$ identifies $V_0$ with the canonical reduced subvariety of 
an irreducible component of $M_0$.  In particular, if $M_0$ is reduced, then $\varphi_0$
is an isomorphism onto an irreducible component.
\end{corollary}

\begin{proof}
This follows from Proposition \ref{first} and Theorem \ref{affine}.
\end{proof}

\vspace{-\baselineskip}
\begin{remark}
When $X$ is a Fano surface and $L=K_X$, this example has an interpretation in string theory.  The quiver variety $M_0$ is the moduli space of vacua for ground states of open strings ending on a D-brane at the tip of the cone $V_0$ (see for example \cite{BP}). Considerations from topological string theory imply that one component of this moduli space should correspond to deformations of the D-brane away from the tip, 
and this component is the one picked out by $V_0$.
\end{remark}

\begin{remark}\label{partial}
Suppose that a character $\chi$ is good for the collection $E_1,\ldots,E_n$
on $X$, in the sense of Section \ref{stable}.  Then $\chi$ is also good for
the collection $\pi^*E_1,\ldots,\pi^*E_n$ on $V$.  (This is because
the quiver for the latter collection is obtained by adding arrows to the quiver
for the original collection, thus making it easier for representations to be semistable.)
The quiver variety
$M_\chi$ for the collection on $V$ is projective over $M_0$, and
the component into which $V$ maps by the tautological map is a partial
resolution of $V_0\subs M_0$.  It's easy to check that this partial resolution
is an isomorphism away from the tip of the cone, and that the fiber over the tip is
isomorphic to the variety $X_\chi$ introduced in Section \ref{stable}.
\end{remark}

\begin{example}\label{dp0}
Let $X = \P^1\times\P^1$ with the exceptional collection $E_1 = \O, E_2 = \O(0,1),
E_3 = \O(1,0)$, and $E_4 = \O(1,1)$.
This collection generates a simple spiral ({\em not} a simple helix) with respect
to the canonical bundle $L = \O(-2,-2)$, and $V_0$ is isomorphic to the quotient
of the conifold $\{xy-zw=0\}\subs\C^4$ by the diagonal action of $\Z/2$.  For a more detailed
exposition of this example, see \cite[\S 4]{BP}.
\end{example}
\end{section}

\begin{section}{A canonical projective quotient}\label{resolutions}
Suppose that $E_1,\ldots,E_n$ are line bundles that generate a simple {\em helix} on
a smooth projective variety $X$.  In other words, we are in the situation of Section \ref{zero}
with $L = K_X$ the canonical bundle and $n = \dim X + 1$.
Let $\Q$ be the Bondal  quiver % with relations
associated to the collection $E_1,\ldots,E_n$ on $X$, and $\Q'$
the Bondal quiver % with relations 
associated to the collection $\pi^*E_1,\ldots,\pi^*E_n$ on 
the total space $V$ of $K_X$.
Then the underlying quiver $Q$ has arrows from $i+1$ to $i$ given by a basis for the vector space
$\Hom_X(E_i,E_{i+1})$, and no arrows between nonadjacent vertices \cite[\S 4]{Br}.
Similarly, Bridgeland shows that
$Q'$ is obtained from $Q$ by adding arrows from $1$ to $n$ given by a basis for
$\Hom_X(E_n,E_1\otimes K_X^{-1})$. 

By Theorem \ref{equiv} and Example \ref{bridge}, the derived functors 
$$RF:\Db\Coh(X)\to\Db\Rep(\Q)\hspace{15pt}\text{and}\hspace{15pt}
 RF':\Db\Coh(V)\to\Db\Rep(\Q')$$
are both equivalences of categories.
Then by Theorem \ref{bz} and the fact that $X$ and $V$ are smooth, 
the tautological maps
$$T:X\to \roq\hspace{15pt}\text{and}\hspace{15pt}
T':V\to \roqp$$ are open immersions,
and therefore the loci of points in $R_{\a}(\Q)$ and $R_{\a}(\Q')$ lying over the images of
these maps are open.
Let $R_{\a}(\Q)_{taut}$ and $R_{\a}(\Q')_{taut}$ be the closures of these loci;
since $X$ is irreducible, $R_{\a}(\Q)_{taut}$ and $R_{\a}(\Q')_{taut}$ are irreducible components
of $R_{\a}(\Q)$ and $R_{\a}(\Q')$.
We will introduce resolutions $\res$ and $\resp$ of 
$R_{\a}(\Q)_{taut}$ and $R_{\a}(\Q')_{taut}$, respectively.  
We will then show that, under certain hypotheses, $\resp$ has the structure
of a $G$-equivariant line bundle over $\res$, and that the GIT
quotient of $\res$ with respect to this line bundle is equal to $X$
in its anticanonical projective embedding. 
We thus recover $X$
as a GIT quotient of a smooth variety by $G$,
without having to make any choice of character.

Since $\a_i=1$ for all $i$, the affine variety $R_\a(\Q)$ admits a particularly simple explicit
description.  We have
\begin{equation*}\label{rq}
R_{\a}(\Q) = \Big\{(a_{ij})\in\Pd_{i,j}\Hom(E_i,E_j)^\vee\bigmid
\ss{\psi_{ijk}^\vee(a_{ik}) = a_{ij}\otimes a_{jk}\hs\text{ for all $i,j,k$}
\\\text{and}\quad a_{ii}\(\operatorname{id}_{E_i}\) = 1}
\Big\},
\end{equation*}
where
$$\psi_{ijk}:\Hom_X(E_i,E_j)\otimes\Hom_X(E_j,E_k)\to\Hom_X(E_i,E_k)$$
is the natural composition map and
$$\psi_{ijk}^\vee:\Hom_X(E_i,E_k)^\vee\to\Hom_X(E_i,E_j)^\vee\otimes\Hom_X(E_j,E_j)^\vee$$
is its dual.
Let
$$\res = \Big\{(a,\ell) \bigmid 
a_{ij}\in\ell_{ij} % \in\P\big(\Hom_X(E_i,E_j)\big)
\text{ and }\Psi_{ijk}^\vee(\ell_{ik}) = \ell_{ij}\otimes\ell_{jk}\Big\}
\subs R_{\a}(\Q)\times
\prod_{i<j}\P\big(\Hom_X(E_i,E_j)\big),
$$
where
$$\Psi_{ijk}^\vee:\P\big(\Hom_X(E_i,E_k)\big)\to \P\big(\Hom_X(E_i,E_j)\big)
\otimes\P\big(\Hom_X(E_j,E_k)\big)$$
is the projectivization of
$\psi_{ijk}^\vee$.
Note that for $\Psi_{ijk}^\vee$ to be well defined we need
$\psi_{ijk}^\vee$ to be injective, or equivalently $\psi_{ijk}$
to be surjective.  This, however, is guaranteed by the fact that $P(\Q)$
is generated by arrows between adjacent nodes.
Note that an
element of $R_\a(\Q)$ is determined by the coordinates
$a_{i\, i+1}$ for all $i<n$, and an element of $\res$ is determined by these data along with
the lines $\ell_{i\, i+1}$, but for notational purposes it is still useful to keep track of $a_{ij}$
and $\ell_{ij}$ for all $i<j$.

The space $\resp$ will be defined in a similar manner, but the
fact that $\Q'$ has loops makes the definition slightly more delicate.
Recall that, for all $i,j$, we have
$$\Hom_V(\pi^*E_i,\pi^*E_j) \cong
\bigoplus_{r\geq 0}\Hom_X\!\(E_i,E_j\otimes K_X^{-r}\).$$
An element $a\in R_{\a}(\Q')$ may be regarded as a collection
$(a_{ij}^r)$, $a_{ij}^r\in\Hom_X\!\(E_i,E_j\otimes K_X^{-r}\)^\vee$, which
satisfies the equations
$$(\psi_{ijk}^{rs})^\vee(a_{ik}^{r+s}) =a_{ij}^r\otimes a_{jk}^s,$$
where $\psi_{ijk}^{rs}$ is the restriction of $\psi_{ijk}$
to the $(r,s)$ graded piece of the product 
$$\Hom_V(\pi^*E_i,\pi^*E_j)
\times\Hom_V(\pi^*E_j,\pi^*E_k).$$
We then define
$$\resp = \Big\{(a,\ell)\bigmid 
a_{ij}^r\in\ell_{ij}^r % \in\P\big(\Hom_X(E_i,E_j\otimes K_X^{-r})\big)
\text{ and }(\Psi_{ijk}^{rs})^\vee(\ell_{ik}^{r+s}) 
= \ell_{ij}^r\otimes\ell_{jk}^s\Big\}
\subs R_{\a}(\Q')_{taut}\times
\prod_{i,j,r}\P\big(\Hom_X(E_i,E_j\otimes K_X^{-r})\big),
$$
where 
$(\Psi_{ijk}^{rs})^\vee$
is the projectivization of
$(\psi_{ijk}^{rs})^\vee$.
Once again, these maps are well defined because
the maps $\psi_{ijk}^{rs}$ are surjective, which follows from
Bridgeland's description of $\Q'$.
As in the case of $\Q$, an
element of $\resp$ is completely determined by the data
\begin{eqnarray*}
a_{i\,i+1}^0\in\ell_{i\,i+1}^0 &\subs&% \Arr_{\Q'}(i+1,i)=
\Hom_V(\pi^*E_i,\pi^*E_{i+1})^\vee_0 \cong\Hom_X(E_i,E_{i+1})^\vee
\hspace{15pt}\text{for all $i<n$},\\
\text{and}\hspace{10pt}
a_{n1}^1\in\ell_{n1}^1 &\subs&% \Arr_{\Q'}(1,n)^\vee=
\Hom_V(\pi^*E_n,\pi^*E_1)^\vee_1 \cong 
\Hom_X(E_n,E_1\otimes K_X^{-1})^\vee,
\end{eqnarray*}
subject to certain relations.

Consider the $G$-equivariant projection from $\resp$ to $\res$ given by remembering
only the degree zero parts $a^0$
and $\ell^0$.

\begin{proposition}\label{line}
Suppose that there exists non-negative integers $\{m_{ij}\}$
such that $$\shom_X(E_n,E_1\otimes K_X^{-1}) \cong
\bigotimes_{i,j}\shom_X(E_i,E_j)^{\otimes m_{ij}}.$$
Then the projection from $\resp$ to $\res$ has the structure of an equivariant line bundle.
\end{proposition}

\begin{proof}
Given an element $(a^0,\ell^0)\in\res$, we will show that 
$\ell_{n1}^1\subs\P\big(\Hom_X(E_n,E_1\otimes K_X^{-1})\big)$
is uniquely determined, and that any point $a_{n1}^1\in\ell_{n1}^1$
extends $(a^0,\ell^0)$ to an element of $\resp$.
To see that $\ell_{n1}^1$ is uniquely determined, let $m = \max_{i,j}\{m_{ij}\}$.
By composing maps of the form $\Psi_{i\,i+1\,i+2}$ as we wrap $m$ times
around the quiver $\Q'$, we obtain an equation
$$\ell_{nn}^m\mapsto\ell_{12}^{\otimes m}\otimes\ldots
\otimes\ell_{n-1\,n}^{\otimes m}
\otimes(\ell_{n1}^1)^{\otimes m}.$$
The right hand side of the above line contains 
$\ell_{12}^{\otimes m_{12}}\otimes\ldots\otimes\ell_{n-1\,n}^{\otimes m_{n-1\,n}}
\otimes\ell_{n1}^1$
as a factor.  The lines $\ell_{12}^{\otimes m_{12}}
\otimes\ldots\otimes\ell_{n-1\,n}^{\otimes m_{n-1\,n}}$
and $\ell_{n1}^1$ both lie in $\P\!\(\Hom_X(E_n,E_1\otimes K_X^{-1}\)$, and the symmetries
of the compositions of the maps $\Psi_{ijk}$ imply that anything in their image
is invariant under the interchanging of these factors.  Thus $\ell_{n1}^1$
must be equal to $\ell_{12}^{m_{12}}\otimes\ldots\otimes\ell_{n-1\,n}^{m_{n-1\,n}}$.

The defining equations for $R_{\a}(\Q')$ are linear in $a_{n1}^0$, and therefore
to see that they are satisfied
by every element $a_{n1}^1\in\ell_{n1}^1$, it will suffice to find
{\em one} element of $\resp$ lying over $(a^0,\ell^0)$ in which $a_{n1}^1$
is nonzero.
Suppose that the image of $a^0$ in the stack $\roq$ is equal
to $T(p)$ for some $p\in X$.  In other words, $a^0$ is obtained
from $p$ by choosing an isomorphism $E_i^\vee\!\!\mid_p\cong \k$
for each $i$.
Let $q\in V$ be a nonzero element of the fiber of $K_X$ at $p$.
Our choices of trivializations of the vector spaces
$E_i^\vee\!\!\mid_p$ induce trivializations of the pullbacks
$\pi^*E_i^\vee\!\!\mid_q$, and thus $T'(q)\in \roqp$ lifts
naturally to an element of $R_{\a}(\Q')$ extending $a^0$.
In Remark \ref{gbgs}, we observed that the helix condition
ensures that $\shom_X(E_i,E_{i+1})$ is generated
by global sections for all $i$.  This observation applies equally well when $i=n$,
so $\shom_X(E_n,E_1\otimes K_X^{-1})$ is generated by global sections, as well.
It follows that $T'(q)$ lifts further to a unique element of $\resp$
extending $(a^0,\ell^0)\in\res$, and that $a_{n1}^1\neq 0$.

We have now shown that the fibers of the map from $\resp$
to $\res$ are vector spaces of dimension at most one, and that the dimension
is equal to one over those elements $(a^0,\ell^0)$ such that $a^0$
lies over the image of $T$.  But this is a dense open condition,
and the dimension of the fiber of an algebraic map is an upper
semicontinuous function.  Hence the dimension of the fiber must always be exactly one.
\end{proof}

\vspace{-\baselineskip}
\begin{example}
Consider the data of Example \ref{p2}, in which $X=\P^2$
and $R_{\a}(\Q)$ is the variety of pairs of parallel vectors in the three dimensional
vector space $\Gamma\(\O(1)\)^\vee$.
The quiver $\Q'$ is obtained by adding arrows from node $1$ to node $3$
indexed by the vector space $\Hom\(\O(2),\O(3)\)\cong\Gamma\(\O(1)\)$,
and the relations tell us that an element of $R_{\a}(\Q')$ is a {\em triple} of
parallel vectors in $\Gamma\(\O(1)\)^\vee$.  The projection from
$R_{\a}(\Q')$ to $R_{\a}(\Q)$ has fibers which are generically lines, but the fiber over
zero is a vector space of dimension $3$.
After blowing up the origin of
$\Gamma\(\O(1)\)^\vee$, 
we obtain $$\res = \Big\{a_{12},a_{23}\in\ell\subs\Gamma\(\O(1)\)^\vee\Big\}
\hspace{10pt}\text{and}\hspace{10pt}
\resp = \Big\{a_{12}^0,a_{23}^0,a_{31}^1\in\ell\subs\Gamma\(\O(1)\)^\vee\Big\},$$
and the fibers of the induced projection are all lines.
\end{example}

If the conclusion of Proposition \ref{line} holds,
let $\mathcal{L}$ be the {\em dual} of the corresponding line bundle.

\begin{theorem}\label{smoothgit}
The GIT quotient of $\res$ with respect to the polarization
$\mathcal{L}$ is isomorphic to $X$, with line bundle $K_X^{-1}$.
\end{theorem}

\begin{proof}
Let $\chi_m$ be the character associated to the collection
$\{m_{ij}\}$ as in Section \ref{stable}, and let 
$\chi = \chi_m + (-1,0,\ldots,0,1)$,
so that $\Echi\cong K_X^{-1}$.  
The fact that $\shom_X(E_i, E_{i+1})$ is generated by global sections ensures that
the natural projection from $\res$ to $R_{\a}(\Q)_{taut}$ is
an isomorphism over $R_{\a}(\Q)^{\chi-st}_{taut}$,
and this isomorphism identifies $\mathcal{L}$
with the restriction of the trivial line bundle on $R_{\a}(\Q)$ twisted by the character $\chi$.
Thus $\res\mod_\mathcal{L}G$ is isomorphic to the tautological component
of $R_{\a}(\Q)\mod_\chi G=\Mchi$, which is isomorphic to $X$ in its anticanonical embedding
by Theorem \ref{gn}, Proposition \ref{great}, and Remark \ref{gbgs}.
\end{proof}

\vspace{-\baselineskip}
\begin{remark}
We note that the collection $\{m_{ij}\}$ satisfying the hypotheses
of Proposition \ref{line} may not be unique.  In Example \ref{p2}, we may take
$m_{12} =1$ and $m_{23}=0$, or vice-versa; these choices give us two different
characters $\chi$ which define the same stability condition on $R_{\a}(\Q)$.
Theorem \ref{smoothgit}, on the other hand, requires no choices.
\end{remark}

\begin{remark}
Since $\resp$ is projective over $R_{\a}(\Q')_{taut}$, they have the same ring
of global functions.  This tells us that
$$\Spec\Gamma(\O^G_{\!\resp}) = \Spec\Gamma(\O^G_{\! R_{\a}(\Q')_{taut}})
= R_{\a}(\Q')_{taut}\mod_{\,0} G,$$
the underlying reduced variety of which isomorphic to $V_0$ by
Theorem \ref{affine}.
Recall that $V_0$ is obtained from $V$ by collapsing the
zero section of $K_X$ to a point.
The GIT quotient $\res\mod_\mathcal{L}G$ is isomorphic to $\Proj\O^G_{\!\resp}$,
which is obtained from $\Spec\Gamma(\O^G_{\!\resp})\cong V_0$ by throwing away the tip of the cone
and dividing by the natural action of $\gm$.  This tells us that any nonreduced structure
of $R_{\a}(\Q')_{taut}\mod_{\,0} G$ must be concentrated at the tip of the cone.  If we had known this
fact {\em a priori}, then it would have constituted an alternative proof of Theorem \ref{smoothgit}.
\end{remark}

% \begin{remark}\label{smooth}
% We have observed that $\res$ is isomorphic to $R_{\a}(\Q)_{taut}$ over the image of
% $T$, and thus $T$ induces a natural inclusion of $X$ into $\res$.  It is
% not hard to check that $\res$ is isomorphic to the total space of a vector bundle
% over $X$, with projection given by taking $(a,\ell)$ to $\ell$, and the above
% inclusion may be thought of as the zero section of this vector bundle.
% The action of $G$ on $\res$ leaves the zero section fixed and scales
% the fibers of this vector bundle, so it is not surprising that we can find
% a polarization with respect to which the GIT quotient is isomorphic to $X$.
% This perspective also shows us that $\res$ is smooth.
% \end{remark}

\begin{remark}
We have used the hypothesis that the collection $E_1,\ldots,E_n$ generates
a simple helix throughout this section, but
we remark that in some examples, our methods may be applied to a simple spiral
as well.  Consider, for example, the collection $\O,\O(1,0),\O(1,1),\O(2,1)$ on
$\P^1\times\P^1$, extended to a simple spiral by the canonical bundle $\O(-2,-2)$.
(Note that this is not the same collection that we used in Example \ref{dp0}.)
In this case $Q$ has arrows only between adjacent nodes with all composition maps surjective, 
and $Q'$ is obtained from $Q$ by adding arrows from the first node to the last.
Thus we may define
$\res$ and $\resp$ exactly as we do in the helix case, and Theorem \ref{smoothgit} will still hold. 
\end{remark}
\end{section}

\footnotesize{

}
\end{document}